# A LARGE DEVIATION INEQUALITY FOR VECTOR FUNCTIONS ON FINITE REVERSIBLE MARKOV CHAINS

By Vladislav Kargin

*Courant Institute of Mathematical Sciences*

Let $S_N$ be the sum of vector-valued functions defined on a finite Markov chain. An analogue of the Bernstein–Hoeffding inequality is derived for the probability of large deviations of $S_N$ and relates the probability to the spectral gap of the Markov chain. Examples suggest that this inequality is better than alternative inequalities if the chain has a sufficiently large spectral gap and the function is high-dimensional.

**1. Introduction.** Suppose that a system evolves according to a Markov chain and that properties of the system are described by a vector-valued function $f$. After a sufficiently long time, the average of the realized values of $f$ converges to its expected value. In many practical situations, it is of great interest to determine how long it takes for the average to converge within specified bounds. In other words, we are interested in estimating the probability of a large deviation of the average from its expected value. Large deviation theory gives the asymptotic rate of convergence but is silent about explicit bounds. In the case of a scalar function, the first explicit estimate of the probability of a large deviation was given by Gillman [9] and was later improved by Dinwoodie [6] and Lezaud [14]. For vector-valued functions, we could proceed by applying one-dimensional estimates to each component of the function. If $S_N$ is a vector with $m$ components and we want to estimate $\Pr\{|S_N| \geq \varepsilon N\}$, then it is enough to estimate $\Pr\{|S_N^i| \geq \varepsilon/\sqrt{m}N\}$, where $S_N^i$ is the $i$th component of the vector sum $S_N$. Since one-dimensional inequalities have the form $\Pr\{|S_N^i| \geq \eta N\} \leq C\exp(-\alpha\eta^2 N)$, our estimate will be

$$\Pr\{|S_N| \geq \varepsilon N\} \leq Cm\exp(-(\alpha/m)\varepsilon^2 N),$$

which has an exponential rate inversely related to $m$. It turns out that it is possible to improve on this inequality by deriving a genuine multidimensional inequality in which the rate function is dimension-free.

To fix notation, let $\mathbb{S}$ be the state space of a finite Markov chain with transition matrix $P$ and invariant distribution $\mu$. We will assume that the chain is *reversible*, that is, that $\mu_s P_{st} = \mu_t P_{ts}$ for any $s$ and $t$ from $\mathbb{S}$. The transition matrix of a reversible chain is similar to a symmetric matrix (i.e., there exists a $D$ such that $D^{-1}PD$ is symmetric) and therefore enjoys many good properties of symmetric matrices. In particular, its eigenvalues are real. Let us denote the eigenvalues of $P$ as $\lambda_i$, where

$$\lambda_0 = 1 > \lambda_1 \geq \lambda_2 \geq \cdots \geq \lambda_{|S|-1} \geq -1.$$

The difference $1 - \lambda_1$ is called the *spectral gap* of the chain. In our study, it will be the main indicator of how well the chain mixes the states. Finally, let $f$ be a function on $\mathbb{S}$ that takes values in an $m$-dimensional real Euclidean space, that is, in a vector space endowed with a scalar product $\langle \cdot, \cdot \rangle$ and the corresponding norm $|\cdot|$. We study the behavior of partial sums $S_N = \sum_{t=1}^{N} f(s_t)$, where the sequence $s_1, \ldots, s_N$ is a realization of the Markov chain evolution.

The behavior of the sum depends on the interaction of properties of the function and the Markov chain. We will use two parameters that characterize this interaction. We call them the *$l^\infty$-norm* and the *principal variance* of $f$. The $l^\infty$-norm is defined as $\|f\|_\infty =: \sup_s |f(s)|$. The principal variance is defined as follows. With each vector $u$, we can associate the variance of the random scalar product $\langle f(s), u\rangle$. The randomness comes from $s$, which is drawn according to the invariant distribution. The *principal variance* of $f$ is defined as the supremum of these variances over all unit vectors $u$:

$$\begin{aligned}\sigma^2(f) &=: \sup_{|u|=1} \sum_{s\in\mathbb{S}} \mu_s \langle f(s), u\rangle^2 \\ &= \sup_{|u|=1} \mathrm{E}[\langle f(s), u\rangle^2].\end{aligned}$$

(In what follows, we will always use symbols E and $\mathrm{E}^{(0)}$ to denote the expectation values relative to the invariant and initial distributions on $\mathbb{S}$, resp.) The principal variance measures the variation of the function $f$ in the long run, when the distribution of $f(s)$ is approximately invariant. The $l^\infty$-norm helps us to determine if the function has an outlier. Directly from the definitions, it is clear that $\sigma^2(f) \leq \|f\|_\infty^2$.

The behavior of the partial sums $S_N$ also depends on the initial distribution $\mu^{(0)}$. It is convenient to use the following measure of the distance between the initial and the invariant distribution:

$$\|\mu^{(0)}/\mu\|^2 =: \mathrm{E}\left[\left(\frac{\mu^{(0)}(s)}{\mu(s)}\right)^2\right].$$

Here is the main result.

THEOREM 1. *Suppose* (1) $P$ *is reversible with spectral gap* $g$, (2) $Ef = 0$, (3) $\|f\|_\infty \le L$ *and* (4) $\sigma^2(f) \le \sigma^2$. *For arbitrary* $\varepsilon > 0$,

$$\Pr\{|S_N| \ge \varepsilon N\} \le 3\|\mu^{(0)}/\mu\| 2^{m/2} e^{-(1/(8k))\varepsilon^2 N}, \tag{1}$$

*where*

$$k = \sigma^2\left(\frac{1}{2} + \frac{1}{g}\right) + L^2 \frac{192}{125} \frac{g}{\log^2[1 + g/2]}.$$

In view of the inequality $\sigma^2(f) \le L^2$, we can take $\sigma^2 = L^2$ and obtain the following estimate that involves only $L$.

COROLLARY 2. *Under the assumptions of Theorem 1,*

$$\Pr\{|S_N| \ge \varepsilon N\} \le 3\|\mu^{(0)}/\mu\| 2^{m/2} \exp\left[-\alpha \frac{\varepsilon^2}{L^2} N\right],$$

*where*

$$\alpha = \left(4 + \frac{8}{g} + \frac{1536}{125} \frac{g}{\log^2[1 + g/2]}\right)^{-1}.$$

REMARKS. 1. Recall that one form of the Bernstein–Hoeffding inequality for i.i.d. and one-dimensional variables is

$$\Pr\{|S_N| \ge \varepsilon N\} \le 2 \exp\left[-\frac{1}{2} \frac{\varepsilon^2}{L^2} N\right] \tag{2}$$

(see, e.g., [10], Theorem 2). This inequality has the same form as the inequality we formulated in Corollary 2, but a better exponential rate. For Markov chains and one-dimensional functions $f$, Gillman [9] showed that if $\|f\|_\infty \le 1$, then

$$\Pr\{S_N \ge \varepsilon N\} \le 2\|\mu^{(0)}/\mu\| \exp\left[-\frac{g}{20\nu} \varepsilon^2 N\right], \tag{3}$$

where $\nu$ is the spread of $P$, that is, $\nu = \max(\mu)/\min(\mu)$. The inequality in Theorem 1 generalizes (3) to the case of multidimensional functions $f$.

2. For a fixed $m$, the probability of large deviations declines exponentially with rate at least $-(8k)^{-1}\varepsilon^2$. Note that this bound on the rate does not depend on the dimension of the Euclidean space where $f$ takes its values. However, the dimension can significantly affect the term before the exponential, which grows exponentially in $m$.

TABLE 1
*Sample size needed to ensure that* $\Pr\{|S_N/N| > 0.01\} < 5\%$

| | Complete graph | | Hypercube | | Circle | |
|---|---|---|---|---|---|---|
| **Method** | $m = 1$ | $m = 20$ | $m = 1$ | $m = 20$ | $m = 1$ | $m = 20$ |
| Theorem 1 | 4 mln | 9 mln | 9 mln | 22 mln | 560 mln | 960 mln |
| Martingale inequality[1] | 280 mln | — | 280 mln | — | 300 mln | — |
| Gillman | 0.7 mln | 26 mln | 2 mln | 80 mln | 160 mln | 2,640 mln |

[1]While [11] derive bounds for vector-valued martingales, they do not provide explicit constants for their inequalities.

EXAMPLES. In the following examples, we study random walks on graphs. We will assume that $E(f) = 0$ and $L = \sigma^2 = 1$. We ask how large $N$ should be to ensure that the following inequality holds:

$$\Pr\{|S_N/N| \geq 0.01\} \leq 0.05.$$

We will consider three examples: a complete graph, a hypercube and a circle. We will set the number of vertices equal to 32 in all examples to make them comparable. (In the example with the circle, we use 33 vertices to ensure that the chain is aperiodic.) We will also assume that the random walks start from the uniform distribution. The results are collected in Table 1.

EXAMPLE 3. *Random walk on a complete graph.* The most connected of all graphs is the complete graph, where each vertex is connected with each of the other vertices. We consider a random walk on a complete graph with $n = 32$ vertices. The spectral gap for this random walk is $n/(n-1) = 1 + 1/31$ (see [1] for derivation).

EXAMPLE 4. *Random walk on a hypercube.* Let the state space be the set of vertices of a 5-dimensional hypercube. With probability 5/6, the next state will be one of the 5 adjacent vertices and with probability 1/6, it remains the same. The spectral gap is $g = 2/(5+1) = 1/3$ (see [5] or [19]).

EXAMPLE 5. *Random walk on a circle.* We also consider a random walk on a circle that consists of $n = 33$ states. If the current state is $x \in \{1, \dots, n\}$, then the next state is $x \pm 1 \bmod (n)$, with probability 1/2 on each possibility. The spectral gap is $g = 1 - \cos(\pi/n) \approx 0.0045$ (see [5] or [19]).

We consider two dimensions, $m = 1$ and $m = 20$, and three methods. The first is from our Theorem 1, the second is given by Gillman's inequality, modified to make it applicable to multidimensional situations, and the third is the method of reduction to martingale inequalities. The following is a sketch of the third method in its application to a random walk on an $n$-vertex graph. Assume that the walk has been started from the invariant distribution. We can define $F_k = E(S_N|s_1, \dots, s_k)$, that is, the expectation of the sum $S_N$ conditional on the first $k$ realizations of the chain. Then $F_1, \dots, F_N$ form a martingale and $F_N = S_N$. For the application of the Bernstein inequality for martingale sequences, we need an estimate on $|F_k - F_{k-1}|$. Using coupling arguments, it is possible to show that $|F_k - F_{k-1}|$ is less than $2(n-1)L$, where $n$ is the number of vertices in the graph. Therefore, for $m = 1$ we have the Bernstein inequality

$$\Pr\{|S_N| \geq \varepsilon N\} \leq 2 \exp\left(-\frac{1}{2} \frac{\varepsilon^2}{[2(n-1)L]^2}\right) \tag{4}$$

and for $m > 1$, similar inequalities are given by Kallenberg and Sztencel [11] (without explicit constants). Note that this method ignores how well the chain mixes and uses only the size of the graph to bound the probability of a large deviation.

Table 1 shows that Gillman's inequality provides the best bounds for $m = 1$, but performs worse than the bound in Theorem 1 for $m = 20$. The martingale inequality underperforms other methods for both the complete graph and hypercube, but is better than the bound in Theorem 1 for the case of the circle. This leads us to the conclusion that the bound in Theorem 1 is most effective for large dimensions and well-connected graphs for which the spectral gap is large.

To put the problem in perspective, we shall sketch a history of the question. Apparently, the first version of a large deviation inequality for sums of i.i.d. random variables was proved by Bernstein in 1924 (see Paper 5 in [3]). Later, Bernstein's result was significantly clarified and improved by Kolmogoroff [13], Chernoff [4], Prokhorov [17], Bennett [2] and Hoeffding [10]. In addition, Hoeffding [10] showed how the inequality can be extended to some classes of dependent variables and, in particular, to the case of martingale differences. Prokhorov [18] proved the multidimensional analogue of the Bernstein inequality for i.i.d. random variables. The multidimensional analogue was also derived by Yurinskii [20] by a different method which is applicable to the case of random variables that take values in an infinite-dimensional Banach space. Later, the multidimensional large-deviation inequalities were generalized to the case of martingale sequences in [11]. They showed that a martingale process with values in a Hilbert space can be represented by a martingale process that takes values in the plane $\mathbb{R}^2$. This device allows reduction of the question of large deviations in many dimensions to the question of large deviations for two-dimensional martingale processes.

For functions defined on the state-space of a finite Markov chain, large deviations were first studied by Miller [15]. Very definitive and general results in this direction were later obtained by Donsker and Varadhan [7]. They established the existence of the exponential rate of the decline in the probability of large deviations and showed how to compute this rate. Their results are valid for vector-valued or even measure-valued functionals of Markov chains acting on very general state spaces. While results of this type are very useful for understanding the asymptotic behavior of large deviations, they do not provide explicit bounds on the probability of a large deviation in a finite sample.

The first one-dimensional Bernstein-type inequality for finite Markov chains was proved by Gillman [9] (see also [6] and [14] for significant improvements). Gillman's method is to write

$$\begin{aligned}\Pr\{S_N \geq \varepsilon N\} &\leq \mathrm{E}^{(0)} \exp(-\theta\varepsilon N + \theta S_N)\\ &= \exp(-\theta\varepsilon N) \sum_{s_0, s_1, \ldots, s_N} \mu^{(0)}_{s_0} P_{s_0 s_1} e^{\theta f(s_1)} \cdots P_{s_{N-1} s_N} e^{\theta f(s_N)}\\ &= \exp(-\theta\varepsilon N)(\mu^{(0)}, [P(\theta)]^N \mathbf{1}_{\mathbb{S}}),\end{aligned}$$

where $P_{st}$ denotes the transition probability from state $s$ to $t$, $P(\theta)$ is a matrix with entries $P_{st} = P_{st} e^{\theta f(t)}$, $\mu^{(0)}$ is the initial distribution, $\mathbf{1}_{\mathbb{S}}$ is a function that takes value 1 on every state of $\mathbb{S}$ and $(\cdot, \cdot)$ denotes a scalar product for functions on $\mathbb{S}$. It turns out that $P(\theta)$ is similar to a symmetric matrix and therefore its norm can be bounded in terms of its eigenvalues. Therefore, the main task is to estimate the eigenvalues of $P(\theta)$, which can be done using Kato's theory of linear operator perturbations. Dinwoodie [6] and Lezaud [14] use a similar method and improve upon Gillman by employing more sophisticated and difficult versions of perturbation theory. Prior to Gillman, the method of a perturbed transition kernel was used by Nagaev [16] to study central limit theorems for Markov chains.

Obviously, Gillman's method is not directly applicable to the case of vector functions since we cannot develop $\mathrm{E}\exp(-\theta\varepsilon N + \theta\|S_N\|)$ in the sum of products of $\exp\|f(s)\|$. To circumvent this difficulty, we use an idea of Prokhorov [18], which was used to prove the multidimensional analogue of the Bernstein inequality for i.i.d. variables. The idea is to consider $\mathrm{E}\exp(-\theta\varepsilon N + \theta\langle S_N, u\rangle)$, where $u$ is a random vector from an appropriate distribution, and later integrate it over the distribution of $u$. The advantage is that $\mathrm{E}\exp(-\theta\varepsilon N + \theta\langle S_N, u\rangle)$ *can* be developed as the sum of products of $\exp\langle f(s), u\rangle$. Using this idea we are able to extend the Bernstein–Gillman inequality to vector functions.

A large body of related literature studies the explicit rates of convergence of a Markov chain to its invariant distribution. For a review, see the book by Diaconis [5], the review paper by Saloff-Coste [19] and the dissertation by Gangolli [8]. Our problem is of a somewhat different flavor because, even for a chain which starts in the invariant distribution, the problem of estimating the probability of a large deviation of the function sum is not trivial.

The rest of the paper is devoted to the proof of the main result. It is organized as follows. Section 2 gives an outline of the proof and explicates the relation of our problem to the eigenvalue problem for a perturbed transition matrix. Section 3 applies a mixture of techniques from the Rellich and Kato perturbation theories to estimate the largest eigenvalue of the perturbed transition matrix. Section 4 concludes.

## 2. Outline of the proof.

Let

$$F_r(x) = \int \exp\langle x, u\rangle \, d\Phi(u),$$

where $x$ and $u$ are vectors from an $m$-dimensional real Euclidean space and $d\Phi(u)$ is the Gaussian measure with density

$$\phi(u) = \frac{1}{(2\pi r^2)^{m/2}} \exp\left(-\frac{|u|^2}{2r^2}\right).$$

We can easily calculate $F_r(x)$ explicitly:

$$F_r(x) = e^{(r^2/2)|x|^2}.$$

Consequently, we can write

$$\begin{aligned} \Pr\{|S_N| \geq \varepsilon N\} &= \Pr\{e^{(r^2/2)|S_N|^2} \geq e^{(r^2/2)|\varepsilon N|^2}\} \\ &\leq e^{(-r^2/2)|\varepsilon N|^2} \mathrm{E}^{(0)}\{e^{(r^2/2)|S_N|^2}\} \\ &= e^{(-r^2/2)|\varepsilon N|^2} \mathrm{E}^{(0)}\left[\int \exp\langle S_N, u\rangle \, d\Phi(u)\right] \\ &= e^{(-r^2/2)|\varepsilon N|^2} \int [\mathrm{E}^{(0)} \exp\langle S_N, u\rangle] \, d\Phi(u). \end{aligned} \tag{5}$$

Consider now $\mathrm{E}^{(0)} \exp\langle S_N, u\rangle$. We will write this expression as a quadratic form and show that what matters is the largest eigenvalue of this form. We will then show that a sufficiently good estimate on the eigenvalue would imply the inequality in Theorem 1. The derivation of the eigenvalue estimate is given in the next section.

Define the *perturbed transition matrix* as a matrix with the following entries:

$$P_{st}(u) = P_{st} e^{\langle f(t), u\rangle}.$$

We denote its largest eigenvalue by $\lambda_0(u)$. Let $(\cdot, \cdot)$ denote the scalar product $(a, b) = \sum_s a_s b_s$, where $s$ are states of the chain and $a_s$ and $b_s$ are scalar-valued functions of $s$. Also, let $\mathbf{1}_{\mathbb{S}}$ denote the scalar-valued function that takes the value 1 on all states.

LEMMA 6.

$$\mathrm{E}^{(0)} \exp\langle S_N, u\rangle = (\mu^{(0)}, [P(u)]^n \mathbf{1}_{\mathbb{S}}).$$

PROOF. We can write

$$\begin{aligned}\mathrm{E}^{(0)} \exp\langle S_N, u\rangle &= \sum_{s_0, s_1, \ldots, s_N} \mu^{(0)}_{s_0} P_{s_0 s_1} e^{\langle f(s_1), u\rangle} \cdots P_{s_{N-1} s_N} e^{\langle f(s_N), u\rangle} \\ &= \sum_{s_0, s_1, \ldots, s_N} \mu^{(0)}_{s_0} P_{s_0 s_1}(u) \cdots P_{s_{N-1} s_N}(u) \\ &= (\mu^{(0)}, [P(u)]^N \mathbf{1}_{\mathbb{S}}). \qquad \square\end{aligned}$$

A fortunate consequence of the reversibility of $P$ is that matrices $P$ and $P(u)$ become symmetric in a coordinate system with dilated axes. This implies that matrices $P$ and $P(u)$ enjoy all of the good properties of symmetric matrices and, in particular, that their eigenvalues are real and their norms can be expressed in terms of the eigenvalue with the largest absolute value.

The second instance of good luck is that both $P$ and $P(u)$ are nonnegative in the sense that all of their entries are nonnegative. This implies that the Perron–Frobenius theorem is applicable and we can pinpoint which of the eigenvalues has the largest absolute value. As we might expect, the largest eigenvalue has the largest absolute value. As a consequence, we are able to estimate the norm of $P(u)$ in terms of its largest eigenvalue and therefore obtain a bound on the value of $(\mu^{(0)}, [P(u)]^N \mathbf{1}_{\mathbb{S}})$.

LEMMA 7. *Let* $D = \mathrm{diag}\{\sqrt{\mu_s}\}$ *and* $E_u = \mathrm{diag}\{\exp \frac{1}{2}\langle f(s), u\rangle\}$. *Define* $S =: DPD^{-1}$ *and* $S_u =: E_u S E_u$. *Then* (1) $S$ *and* $S_u$ *are symmetric,* (2) $S_u$ *is similar to* $P(u)$ *and has the same eigenvalues as* $P(u)$, (3) *the eigenvalues of* $P(u)$ *are real and* (4) *the largest eigenvalue of* $P(u)$ *has the largest absolute value among all eigenvalues of* $P(u)$.

REMARK. Here, $S$ and $E_u$ denote matrices and should not be confused with the notation for the Markov chain, $\mathbb{S}$, and for the expectation value, E, respectively.

PROOF OF LEMMA 7. First, the reversibility of $P$ implies that $S =: DPD^{-1}$ is symmetric. Indeed,

$$\begin{aligned}S_{ji} \equiv \mu_j^{1/2} P_{ji} \mu_i^{-1/2} &= \mu_j^{-1/2} \mu_j P_{ji} \mu_i^{-1/2} \\ &= \mu_j^{-1/2} P_{ij} \mu_i \mu_i^{-1/2} = S_{ij}.\end{aligned}$$

Then $S_u = E_u S E_u$ is symmetric because $E_u$ is symmetric. It is similar to $P(u)$ because

$$\begin{aligned} P(u) &\equiv PE_u^2 = D^{-1}SDE_u^2 \\ &= D^{-1}E_u^{-1}(E_u S E_u)E_u D \\ &= (E_u D)^{-1} S_u (E_u D), \end{aligned}$$

where we have used the commutativity of $D$ and $E_u$. Consequently, $S_u$ and $P(u)$ have the same eigenvalues. The eigenvalues of $S_u$ are real because $S_u$ is symmetric. Therefore, the eigenvalues of $P(u)$ are also real. Finally, $P(u)$ has nonnegative entries and therefore, by the Perron–Frobenius theorem, its largest eigenvalue has the largest absolute value. □

LEMMA 8. *If the chain $P$ is reversible, $|u| \le 1$ and $|f(s)| \le 1$ for any $s$, then*

$$(\mu^{(0)}, [P(u)]^N \mathbf{1}_{\mathbb{S}}) \le 3\|\mu^{(0)}/\mu\| \lambda_0(u)^N.$$

PROOF. Since $S_u$ is symmetric and its largest eigenvalue has the largest absolute value, then $\|S_u\| \le \lambda_0(u)$. Therefore:

$$\begin{aligned} (\mu^{(0)}, [P(u)]^N \mathbf{1}) &= (\mu^{(0)}(E_u D)^{-1}, S_u^N (E_u D)\mathbf{1}_{\mathbb{S}}) \\ &\le \lambda_0(u)^N \|\mu^{(0)}(E_u D)^{-1}\| \|(E_u D)\mathbf{1}_{\mathbb{S}}\|, \end{aligned}$$

where $\|\cdot\|$ denotes the norm corresponding to the scalar product $(\cdot,\cdot)$.

Then

$$\begin{aligned} \|\mu^{(0)}(E_u D)^{-1}\| &= \left(\sum_s \frac{[\mu_s^{(0)}]^2}{\mu_s} \exp\langle -f(s), u\rangle\right)^{1/2} \\ &\le \sqrt{3} \sum_s \frac{[\mu_s^{(0)}]^2}{\mu_s} \\ &= \sqrt{3}\|\mu^{(0)}/\mu\|, \end{aligned}$$

where we have used the fact that $|\langle f(s), u\rangle| \le |f(s)||u| \le 1$ and consequently $\exp\langle \pm f(s), u\rangle \le 3$. Similarly,

$$\|(E_u D)\mathbf{1}_{\mathbb{S}}\| = \left(\sum_s \mu_s \exp\langle f(s), u\rangle\right)^{1/2} \le \sqrt{3}.$$

Combining, we get

$$(\mu^{(0)}, [P(u)]^N \mathbf{1}_{\mathbb{S}}) \le 3\|\mu^{(0)}/\mu\| \lambda_0(u)^N. \qquad \square$$

Suppose, for the moment, that we have managed to establish the inequality

$$\lambda_0(u) \leq \exp(k|u|^2).$$

Then, using Lemmas 6 and 8, we can write

$$\begin{aligned}\int [\mathrm{E}\exp\langle S_N, u\rangle]\, d\Phi(u) &\leq 3\|\mu^{(0)}/\mu\| \int \lambda_0(u)^N \, d\Phi(u)\\ &\leq 3\|\mu^{(0)}/\mu\| \int \exp(k|u|^2 N)\, d\Phi(u)\\ &= \frac{3\|\mu^{(0)}/\mu\|}{(2\pi r^2)^{m/2}} \int \exp(k|u|^2 N) e^{-|u|^2/(2r^2)}\, du.\end{aligned}$$

In spherical coordinates, we can rewrite this expression as follows:

$$\begin{aligned}&\frac{3\|\mu^{(0)}/\mu\|}{(2\pi r^2)^{m/2}} \frac{m\pi^{m/2}}{\Gamma(m/2+1)} \int_0^\infty t^{m-1} e^{kt^2N - t^2/(2r^2)}\, dt\\ &\quad = \frac{3\|\mu^{(0)}/\mu\| m}{2^{m/2}\Gamma(m/2+1)} \int_0^\infty s^{m-1} e^{kr^2s^2N - s^2/2}\, ds,\end{aligned}$$

where we use the fact that the surface area of the unit sphere in $m$-dimensional real Euclidean space is $m\pi^{m/2}/\Gamma((m/2)+1)$. Next, set

$$r = (2\sqrt{kN})^{-1}. \tag{6}$$

Then

$$\int_0^\infty s^{m-1} e^{kr^2s^2N} e^{-s^2/2}\, ds = \int_0^\infty s^{m-1} e^{-s^2/4}\, ds.$$

Making the substitution $t = s^2/4$, we compute

$$\begin{aligned}\int_0^\infty s^{m-1} e^{-s^2/4}\, ds &= 2^{m-1} \int_0^\infty t^{m/2-1} e^{-t}\, dt\\ &= 2^{m-1}\Gamma\left(\frac{m}{2}\right).\end{aligned}$$

So, combining, we obtain

$$\begin{aligned}\int [\mathrm{E}\exp\langle S_N, u\rangle]\, d\Phi(u) &\leq \frac{3\|\mu^{(0)}/\mu\| m 2^{m-1}\Gamma(m/2)}{2^{m/2}\Gamma(m/2+1)}\\ &= 3\|\mu^{(0)}/\mu\| 2^{m/2}.\end{aligned}$$

Substituting this and (6) into (5), we obtain

$$\Pr\{|S_N| \geq \varepsilon N\} \leq 3\|\mu^{(0)}/\mu\| 2^{m/2} e^{-1/(8k)\varepsilon^2 N},$$

which is the desired inequality.

In the above, we have assumed that $\|f\|_\infty \le 1$. In the general case when $\|f\|_\infty \le L$, we simply introduce the auxiliary function $g = f/L$. Then

$$\Pr\{|f_1 + \cdots + f_N| \ge \varepsilon N\} = \Pr\left\{|g_1 + \cdots + g_N| \ge \frac{\varepsilon}{L} N\right\}$$

and the latter probability can be estimated if we observe that $\|g\|_\infty \le 1$ and $\sigma^2(g) = \sigma^2(f)/L^2$.

It remains to derive the required estimate on the eigenvalue $\lambda_0(u)$.

**3. A bound on the largest eigenvalue of the perturbed transition matrix.** We need to estimate the largest eigenvalue of the perturbed transition matrix $P(u) = P\,\mathrm{diag}(\exp\langle f(t), u\rangle)$. In the following, we use the notation $P(z) = P(zu)$, where $u$ is a fixed vector of length 1. Our main concern will be real values of $z$ which lie in the interval $[0, \infty)$, but we will also need to consider the complex values of $z$. It is known that if $\lambda_i$ is an eigenvalue of $P$ of multiplicity 1, then there is a complex-analytic function $\lambda_i(z)$ defined in a neighborhood of $z = 0$ such that $\lambda_i(z)$ is an eigenvalue of $P(z)$. This function is called the *perturbation of the eigenvalue* $\lambda$. We will consider this function for $i = 0$.

It will be clear from the following discussion that for all sufficiently small $z$, say, for $|z| \le r_\Gamma$, there exists a circle around $\lambda_0(z)$ such that $P(z)$ has no eigenvalues in this circle except $\lambda_0(z)$ itself. Since for real positive $z$, the largest eigenvalue of $P(z)$ must be real and positive (by the Perron–Frobenius theorem) and since initially at $z = 0$, $\lambda_0$ is the largest eigenvalue, we can conclude by continuity that when $z$ changes from zero to $r_\Gamma$ along the real line, the largest eigenvalue of $P(z)$ remains $\lambda_0(z)$. Therefore, for this range of $z$, the desired estimate for the largest eigenvalue of $P(z)$ follows from an appropriate estimate for $\lambda_0(z)$. This estimate will be obtained from Kato perturbation theory. For larger values of the perturbation parameter $z$, we will use a different method which bounds all eigenvalues of $P(z)$ at once.

We know that $\lambda_0(0) = 1$ and it is easy to show that $\lambda_0'(0) = 0$. It is also relatively easy to bound the second derivative of $\lambda_0(z)$ at $z = 0$. It is somewhat more difficult to estimate the remainder $\lambda_0(z) - 1 - \lambda_0''(0)z^2$ in an open neighborhood of $z = 0$. We will establish an estimate by studying the resolvent of the perturbed operator in the complex $z$-plane (the Kato method, see [12]).

For convenience, we shall call the following set of conditions Assumption A:

1. $P$ is a reversible chain with spectral gap $g$;
2. $Ef(s) = 0$;

3. The principal variance of $f$ is $\sigma^2$;
4. $|f(s)| \le 1$ for each $s$.

In the following, we always suppose that Assumption A holds. The main result of this section is the following estimate.

PROPOSITION 9.

$$\lambda_0(v) \le e^{k|v|^2},$$

*where*

$$k = \sigma^2\left(\frac{1}{2} + \frac{1}{g}\right) + \frac{192}{125}\frac{g}{\log^2[1+g/2]}. \tag{7}$$

First, we estimate $\lambda_0'(0)$ and $\lambda_0''(0)$.

LEMMA 10. $\lambda_0'(0) = 0$.

PROOF. Matrix $P(z)$ can be developed as a power series in $z$:

$$P(z) = P\left(\sum_{n=0}^{\infty} \frac{1}{n!} V^n z^n\right), \tag{8}$$

where

$$V = \mathrm{diag}\{\langle f(t), u\rangle\}.$$

Let the expansions for $\lambda_0(z)$ and the corresponding eigenvector, $X(z)$, be

$$\lambda_0(z) = 1 + \lambda'(0)z + \tfrac{1}{2}\lambda''(0)z^2 + \cdots,$$
$$X(z) = \mu + X'(0)z + \tfrac{1}{2}X''(0)z^2 + \cdots.$$

Writing the equality $X(z)P(z) = \lambda_0(z)X(z)$ in powers of $z$, we obtain

$$\begin{aligned} \mu P &= \mu, \\ X'(0)P + \mu P V &= \lambda'(0)\mu + X'(0). \end{aligned} \tag{9}$$

Multiply the last line by $\mathbf{1}_{\mathbb{S}}$ on the right and use the facts that $P\mathbf{1}_{\mathbb{S}} = \mathbf{1}_{\mathbb{S}}$ and $\mu\mathbf{1}_{\mathbb{S}} = 1$. (Recall that $\mathbf{1}_{\mathbb{S}}$ is a scalar-valued function that takes the value 1 on all states.) We then obtain

$$\lambda'(0) = \mu V \mathbf{1}_{\mathbb{S}}.$$

However,

$$\mu V \mathbf{1}_{\mathbb{S}} = \sum_s \mu_s \langle f(s), u\rangle \tag{10}$$

$$= (Ef, u) = 0, \tag{11}$$

by assumption. Therefore, $\lambda'(0)=0$. □

We also require some information about the perturbation of the eigenvector, in particular, about $X'(0)$. From (9), $X'(0)$ must satisfy the following equation:

$$X'(0)(I-P)=\mu V. \tag{12}$$

It is tempting to write $X'(0)=(I-P)^{-1}\mu V$. However, $I-P$ is not invertible, which is reflected, for example, in the fact that if a vector $X'$ satisfies equation (12), then $X'+a\mu$ also satisfies it. We need to impose one additional constraint to determine the solution. We choose a normalization in which $X'(0)$ is the unique solution of (12) that satisfies the additional constraint that $(X'(0),\mathbf{1}_{\mathbb{S}})=0$.

To solve (12), we need a pseudo-inverse of $I-P$. The traditional pseudo-inverse is not appropriate because, first, $P$ is not symmetric and second, we use a nonstandard normalization of the solution. An appropriate concept of the pseudo-inverse is as follows.

Let $\mathbf{1}_{\mathbb{S}}^{\perp}$ be the subspace of vectors orthogonal to $\mathbf{1}_{\mathbb{S}}$. This subspace is invariant under the right action of $P$. Indeed, if $x\mathbf{1}_{\mathbb{S}}=0$, then $xP\mathbf{1}_{\mathbb{S}}=x\mathbf{1}_{\mathbb{S}}=0$. We define the pseudo-inverse operator $(I-P)^{\dagger}$ as the inverse of $I-P$ on $\mathbf{1}_{\mathbb{S}}^{\perp}$ and as 0 on $\mathbf{1}_{\mathbb{S}}$.

If $P$ is reversible, then $P=D^{-1}SD$, where $S$ is symmetric. Since the subspace $\mathbf{1}_{\mathbb{S}}^{\perp}$ is invariant under the right action of $P$, the subspace $\mathbf{1}_{\mathbb{S}}^{\perp}D^{-1}$ is invariant under the right action of $S$ and we can define $(I-S)^{\dagger}$, which is the inverse of $I-S$ on $\mathbf{1}_{\mathbb{S}}^{\perp}D^{-1}$ and is zero on $\mathbf{1}_{\mathbb{S}}D^{-1}$. Note that $(I-S)^{\dagger}$ and $S$ commute and that $D^{-1}(I-S)^{\dagger}D=(I-P)^{\dagger}$.

Lemma 11. $X'(0)=\mu V(I-P)^{\dagger}=\mu VD^{-1}(I-S)^{\dagger}D$.

Proof. By (10), $\mu V\in\mathbf{1}_{\mathbb{S}}^{\perp}$. Therefore, the product $\mu V(I-P)^{\dagger}$ satisfies equation (12) and belongs to $\mathbf{1}_{\mathbb{S}}^{\perp}$. Consequently, it coincides with $X'(0)$. □

Now, consider the second derivative of the eigenvalue function.

Lemma 12. $\lambda_0''(0)\le(1+2/g)\sigma^2$.

Proof. Let us equate $z^2$ terms in the expansion of the equality $X(z)P(z)=\lambda(z)X(z)$, taking into account that $\lambda'(0)=0$ and $\mu P=\mu$:

$$\tfrac{1}{2}X''(0)P+X'(0)PV+\tfrac{1}{2}\mu V^2=\tfrac{1}{2}\lambda_0''(0)\mu+\tfrac{1}{2}X''(0).$$

Multiplying this equality by $\mathbf{1}_{\mathbb{S}}$ on the right and using the fact that $P\mathbf{1}_{\mathbb{S}}=\mathbf{1}_{\mathbb{S}}$, we obtain the following formula for $\lambda''(0)$:

$$\lambda_0''(0)=\mu V^2\mathbf{1}_{\mathbb{S}}+2X'(0)PV\mathbf{1}_{\mathbb{S}}. \tag{13}$$

Consider the absolute value of the second term in (13):

$$\begin{aligned}|X'(0)PV\mathbf{1}_{\mathbb{S}}| &= |\mu V D^{-1}(I-S)^{\dagger} SDV\mathbf{1}_{\mathbb{S}}| \\ &\leq \|(I-S)^{\dagger}S\| \|\mu V D^{-1}\| \|DV\mathbf{1}_{\mathbb{S}}\|,\end{aligned}$$

where we used Lemma 11 and the equality $P = D^{-1}SD$. [Here, we use $\|\cdot\|$ to denote both the norm of a function on $\mathbb{S}$ and the norm of an operator that acts on these functions: by definition, $\|f\| = (f,f)^{1/2}$ and $\|A\| = \sup_{\|f\|=1} \|Af\|$.]

The operator $(I-S)^{\dagger}S$ is symmetric with eigenvalues which are either zeros or $\lambda_i/(1-\lambda_i)$, where $i \geq 1$. Consequently,

$$\|(I-S)^{\dagger}S\| \leq \frac{1}{g}.$$

Next,

$$\|DV\mathbf{1}_{\mathbb{S}}\| = \left(\sum_s \mu_s \langle f(s), u\rangle^2\right)^{1/2} \leq \sigma$$

and

$$\|\mu V D^{-1}\| = \left(\sum_s \mu_s \langle f(s), u\rangle^2\right)^{1/2} \leq \sigma,$$

where we used the fact that $D = \operatorname{diag}\{\sqrt{\mu_s}\}$. Combining, we have

$$|X'(0)PV\mathbf{1}_{\mathbb{S}}| \leq \frac{\sigma^2}{g}.$$

Finally, for the first term on the right-hand side of (13), we have

$$|\mu V^2 \mathbf{1}_{\mathbb{S}}| = \left|\sum_s \mu_s \langle f(s), u\rangle^2\right| \leq \sigma^2$$

and therefore

$$\lambda_0''(0) \leq \sigma^2\left(1 + \frac{2}{g}\right). \qquad \square$$

We now turn to the estimation of the residual $\lambda_0(z) - 1 - \lambda_0''(0)z^2$. The following is a quick excursion in Kato's theory of perturbations. The *resolvent* of the perturbed operator $P(z)$ is defined as $R(\zeta, z) \equiv [P(z) - \zeta]^{-1}$. We want to estimate the change in eigenvalues of $P(z)$ when $z$ changes. For this purpose, we study how the resolvent of $P(z)$ depends on $z$.

Let us, for economy of space, write

$$A(z) =: P(z) - P = P(Vz + \tfrac{1}{2}V^2z^2 + \cdots).$$

We can write

$$\begin{aligned} P(z) - \zeta &= P - \zeta + A(z) \\ &= (P - \zeta)[1 + R(\zeta)A(z)] \end{aligned}$$

and consequently,

$$R(\zeta, z) = [1 + R(\zeta)A(z)]^{-1} R(\zeta).$$

The power series for $[1 + R(\zeta)A(z)]^{-1}$ is

$$[1 + R(\zeta)A(z)]^{-1} = \sum_{n=0}^{\infty} [R(\zeta)A(z)]^n. \tag{14}$$

$R(\zeta, z)$ is nonsingular if this power series is convergent, which holds if

$$\|R(\zeta)A(z)\|_{sp} < 1, \tag{15}$$

where $\|\cdot\|_{sp}$ denotes the *spectral norm*,

$$\|X\|_{sp} =: \limsup_{n\to\infty} \|X^n\|^{1/n}.$$

Recall that the reversibility of $P$ implies that it can be represented as $P = D^{-1}SD$, where $D = \mathrm{diag}\{\sqrt{\mu_s}\}$ and $S$ is symmetric. Let us denote $(S - \zeta)^{-1}$ by $R_S(\zeta)$.

LEMMA 13. *The power series* (14) *for* $[1 + R(\zeta)A(z)]^{-1}$ *converges if* $|z| < \log(1 + \|R_S(\zeta)S\|^{-1})$.

PROOF. In our case, the perturbation is

$$A(z) = P(e^{zV} - 1),$$

where $V = \mathrm{diag}(\langle f(s), u\rangle)$. By criterion (15), we should determine when $\|R(\zeta)P(e^{zV} - 1)\|_{sp} < 1$. For reversible $P$, we can write

$$\begin{aligned} R(\zeta)P &\equiv (P - \zeta)^{-1}P \\ &= D^{-1}(S - \zeta)^{-1}SD \\ &= D^{-1}R_S(\zeta)SD. \end{aligned}$$

Using the fact that both $D$ and $(e^{zV} - 1)$ are diagonal and therefore commute, we can further write

$$R(\zeta)P(e^{zV} - 1) = D^{-1}R_S(\zeta)S(e^{zV} - 1)D.$$

Next, we use the property of the spectral norm that it is not changed by similarity transformations and write

$$\begin{aligned} \|R(\zeta)P(e^{zV} - 1)\|_{sp} &= \|R_S(\zeta)S(e^{zV} - 1)\|_{sp} \\ &\le \|R_S(\zeta)S(e^{zV} - 1)\|, \end{aligned}$$

where we also used the fact that the spectral norm is bounded from above by the usual operator norm. We can continue as follows:

$$\|R_S(\zeta)S(e^{zV}-1)\| \leq \|R_S(\zeta)S\| \sum_{k=1}^{\infty} \frac{1}{k!}|z|^k \|V^k\|.$$

From assumptions on $u$ and $f(s)$, it follows that $\|V\| \leq 1$ and consequently,

$$\|R_S(\zeta)S(e^{zV}-1)\| \leq \|R_S(\zeta)S\|(e^{|z|}-1).$$

This expression is less than 1, provided that $|z| < \log(1+\|R_S(\zeta)S\|^{-1})$. □

In the following, it is useful to keep in mind the distinction between the $\zeta$-plane, where the spectral parameter $\zeta$ lives, and the $z$-plane, where the perturbation parameter $z$ lives.

LEMMA 14. *Let $\Gamma$ be a circle of radius $r_\zeta$ in the $\zeta$-plane whose interior contains exactly one eigenvalue of $P$, $\lambda_0 = 1$. Define*

$$r_z = \min_{\zeta \in \Gamma} \log(1+\|R_S(\zeta)S\|^{-1}).$$

*Then for every $z$ in the $z$-plane such that $|z| \leq r_z$, there is exactly one eigenvalue of $P(z)$ inside $\Gamma$ [i.e., the eigenvalue $\lambda_0(z)$ of the perturbed matrix remains inside $\Gamma$].*

*Moreover, for $\alpha \in (0,1)$, the eigenvalue function $\lambda_0(z)$ is holomorphic in the disc $|z| \leq (1-\alpha)r_z$ and its third derivative inside the disc can be estimated as follows:*

$$|\lambda_0'''(z)| \leq \frac{12}{\alpha^3}\frac{r_\zeta}{r_z^3}.$$

Intuitively, if the resolvent $R_S(\zeta)$ is small in magnitude, then we can be sure that for perturbations less then $r_z$, the eigenvalue $\lambda_0(z)$ does not move far from $\lambda_0(0)$ and there are no other eigenvalues near $\lambda_0(z)$. The size of $r_z$ is inversely related to the size of $R_S(\zeta)$.

PROOF OF LEMMA 14. Let $D$ be a circle in the $z$-plane with center at 0 and radius $r_z = \log(1+\|R_S(\zeta)S\|^{-1})$. Consider an arbitrary $z_0$ inside $D$. We can connect $z=0$ and $z_0$ by a curve $\Lambda$ that lies completely inside the circle $D$. When we change $z$ along this curve, the eigenvalues of the operator $P(z)$ follow paths that never intersect the circle $\Gamma$—we know this because by Lemma 13, the power series for the resolvent $R(\zeta, z)$ always converge for all $\zeta \in \Gamma$. Consequently, the number of eigenvalues of the operator $P(z)$ that are located inside $\Gamma$ is conserved along the path $\Lambda$. It follows that $P(z_0)$ has exactly one eigenvalue inside $\Gamma$.

For the second part of the lemma, take an arbitrary $z_0$ such that $|z_0| \leq (1-\alpha)r_z$. Then exactly one eigenvalue of $P(z_0)$ lies inside $\Gamma$. Consider the circle $D_0$ with center at $z_0$ and radius $\alpha r_z$. This circle lies entirely inside the circle $D$ and consequently, for any $z \in D_0$, there is only one eigenvalue of $P(z)$ inside $\Gamma$. Hence,

$$|\lambda_0(z) - \lambda_0(z_0)| \leq 2r_\zeta.$$

Recalling that $\lambda(\varkappa)$ is holomorphic (see [12]), we can estimate its third derivative at $z_0$ by using Cauchy's inequality:

$$|\lambda_0'''(z_0)| \leq 6\frac{\max_{z\in D_0}|\lambda(z)-\lambda(z_0)|}{|z-z_0|^3} = 6\frac{2r_\zeta}{(\alpha r_z)^3} = \frac{12}{\alpha^3}\frac{r_\zeta}{r_z^3}. \qquad \square$$

Lemma 15. *Let $\Gamma$ be a circle of radius $r_\Gamma = g/2$ around $\lambda_0 = 1$. Then*

$$\max_{\zeta\in\Gamma}\|R_S(\zeta)S\| = \frac{2}{g}.$$

Proof. Since $S$ is similar to $P$, it has the same eigenvalues. Since $S$ is symmetric, $R_S(\zeta)S$ is also symmetric and its norm coincides with the largest absolute value of its eigenvalues. Further, $R_S(\zeta)S$ has eigenvalues $(\lambda_i - \zeta)^{-1}\lambda_i$. It is easy to see that if $\zeta \in \Gamma$, then the maximum is reached for $i = 0$ and $\zeta_0 = 1 - g/2$. A calculation gives

$$\|R_S(\zeta_0)S\| = \frac{2}{g}. \qquad \square$$

Lemma 16. *Take $\alpha \in (0,1)$. Then for any $z$ in the disc $|z| \leq (1-\alpha)\log(1+g/2)$, the following inequality holds:*

$$|\lambda_0'''(z)| \leq \frac{6g}{\alpha^3}\log^{-3}\left[1+\frac{g}{2}\right].$$

Proof. From Lemma 14,

$$|\lambda_0'''(z)| \leq \frac{12}{\alpha^3}\frac{r_\zeta}{r_z^3}.$$

Take $r_\zeta = g/2$ and apply Lemma 15 to obtain

$$r_z = \min_{\zeta\in\Gamma}\log(1+\|R_S(\zeta)S\|^{-1}) = \log\left(1+\frac{g}{2}\right).$$

Therefore,

$$|\lambda_0'''(z)| \leq \frac{6g}{\alpha^3}\log^{-3}\left(1+\frac{g}{2}\right). \qquad \square$$

Combining the previous lemmas, we obtain the following result.

LEMMA 17. *Take $\alpha \in (0,1)$. Then for any $z$ in the disc $|z| \leq (1-\alpha)\log(1+g/2)$, the following inequality holds:*

$$|\lambda_0(z)| \leq e^{k|z|^2},$$

*where*

$$k = \sigma^2\left(\frac{1}{2} + \frac{1}{g}\right) + \frac{1-\alpha}{\alpha^3}\frac{g}{\log^2[1+g/2]}.$$

PROOF. First, using Lemmas 12 and 16, we write

$$\begin{aligned}|\lambda_0''(z)| &\leq \sigma^2\left(1 + \frac{2}{g}\right) + \left|\int_0^z \lambda_0'''(t)\,dt\right| \\ &= \sigma^2\left(1 + \frac{2}{g}\right) + \frac{6}{\alpha^3} g \log^{-3}\left[1 + \frac{g}{2}\right]|z|.\end{aligned}$$

Then, using Lemma 10, we get

$$\begin{aligned}|\lambda_0'(z)| &\leq \int_0^{|z|} |\lambda_0''(t)|\,dt \\ &\leq \sigma^2\left(1 + \frac{2}{g}\right)|z| + \frac{3}{\alpha^3} g \log^{-3}\left[1 + \frac{g}{2}\right]|z|^2\end{aligned}$$

and

$$\begin{aligned}|\lambda_0(z)| &\leq 1 + \int_0^{|z|} |\lambda_0'(t)|\,dt \\ &\leq 1 + \sigma^2\left(\frac{1}{2} + \frac{1}{g}\right)|z|^2 + \frac{1}{\alpha^3} g \log^{-3}\left[1 + \frac{g}{2}\right]|z|^3.\end{aligned}$$

Using the condition $|z| \leq (1-\alpha)\log[1+g/2]$, we further reduce this to

$$|\lambda_0(z)| \leq 1 + \left[\sigma^2\left(\frac{1}{2} + \frac{1}{g}\right) + \frac{1-\alpha}{\alpha^3} g \log^{-2}\left[1 + \frac{g}{2}\right]\right]|z|^2.$$

This inequality and the inequality $1 + x^2 \leq e^{x^2}$ together imply the claim of the lemma. □

We should now treat the case when $z$ is real and greater than $(1-\alpha)\log(1+g/2)$.

LEMMA 18. *For every real $z > 0$,*

$$|\lambda_0(z)| \leq e^z.$$

PROOF. Recall (from Lemma 7) that $P(z)$ has the same eigenvalues as $S(z)$, where $S(z) = E_{z/2} S E_{z/2}$, $E_{z/2} = \operatorname{diag} \exp(\frac{z}{2}\langle f(t), u\rangle)$ and $u$ is a vector of unit length. It follows that the absolute value of the largest eigenvalue does not exceed $\|S(z)\| \leq \|E_{z/2}\|^2 \|S\| \leq e^z$, where we used the assumption that $|f(t)| \leq 1$ to bound $\|E_{z/2}\|$. □

LEMMA 19. *For every real* $z > 0$,

$$|\lambda_0(z)| \leq e^{k|z|^2}, \tag{16}$$

*where*

$$k = \sigma^2 \left(\frac{1}{2} + \frac{1}{g}\right) + \frac{192}{125} \frac{g}{\log^2[1 + g/2]}. \tag{17}$$

PROOF. Take $\alpha = 5/8$. Then by Lemma 17, inequality (16) with rate (17) holds for $|z| \leq (3/8) \log(1 + g/2)$. However, for $|z| \geq (3/8) \log(1 + g/2)$, we have

$$\begin{aligned} k|z|^2 &\geq \left(\sigma^2 \left(\frac{1}{2} + \frac{1}{g}\right) + \frac{192}{125} \frac{g}{\log^2[1 + g/2]}\right) \frac{3}{8} \log(1 + g/2)|z| \\ &\geq \frac{72}{125} \frac{g}{\log[1 + g/2]} |z| \geq |z| \end{aligned}$$

and using Lemma 18, we conclude that inequality (16) with rate (17) is valid for all real $z > 0$. □

The claim of Proposition 9 follows if we take $z = |v|$ and $u = v/|v|$ in Lemma 17. As was shown in the previous section, the validity of the inequality in Proposition 9 implies the validity of Theorem 1.

**4. Concluding remarks.** We have derived an inequality for the probability of large deviations of vector-valued functions on a finite Markov chain. The results can be extended in two directions. First, it is desirable to eliminate dependence on the dimension in the term before the exponential. Corresponding results for i.i.d. and martingale variables suggest that this is possible. Second, it would be desirable to extend the results to denumerable Markov chains and, in particular, to random walks on denumerable groups

**Acknowledgment.** I would like to thank Diana Bloom for her editorial help.

Courant Institute
of Mathematical Sciences
109-20 71st Road
Apt. 4A
Forest Hills, New York
USA
E-mail: kargin@cims.nyu.edu